\documentclass[11pt,a4paper,twoside]{scrartcl} 

\usepackage[includemp,
            a4paper,
            inner=2.2cm,
            outer=2.2cm,
            top=2.5cm,
            bottom=3cm,
            marginparwidth=4.5cm,
            marginparsep=0.4cm]{geometry}

\newlength{\overflowingheadlen}
\usepackage{scrlayer-scrpage}
\renewpagestyle{scrheadings}{
  {\hspace{-\marginparwidth}\hspace{-\marginparsep}%
   \makebox[\overflowingheadlen][l]{\makebox[2em][r]{\thepage}\quad\rule{0.5pt}{9pt}\quad{}\leftmark}}%
  {\makebox[\overflowingheadlen][r]{\rightmark\quad\rule{0.5pt}{9pt}\quad\makebox[2em][l]{\thepage}}}%
  {}
}{
  {}{}{}
}
\renewpagestyle{plain.scrheadings}{
  {}{}{}
}{
  {\thepage}%
  {\makebox[\overflowingheadlen][r]{\thepage}}%
  {}
}

\usepackage[most]{tcolorbox}

\newtcolorbox{aside}{
  colback=gray!10,      
  colframe=gray!10,     
  boxrule=0pt,          
  arc=2pt,              
  left=8pt, right=8pt, top=6pt, bottom=6pt,
  before skip=10pt, after skip=10pt
}

\usepackage[utf8]{inputenc}
\usepackage[T1]{fontenc}
\usepackage[english]{babel}
\usepackage[biolinum]{libertine}  
\setkomafont{disposition}{\sffamily\bfseries} 
\setkomafont{author}{\normalsize}
\setkomafont{date}{\normalsize}

\usepackage[hypcap=true]{caption}
\usepackage{marginfix}
\extrafloats{100}

\usepackage{floatrow}
\usepackage[footnote]{snotez}
\setsidenotes{text-mark-format=\textsuperscript{\normalfont#1},
              note-mark-format=#1:,
              note-mark-sep=\enskip}

\usepackage{graphicx}
\graphicspath{{figs/}}
\PassOptionsToPackage{hyphens}{url}
\usepackage[breaklinks=true,colorlinks=true,linkcolor=black,citecolor=black,urlcolor=blue]{hyperref}

\usepackage{xcolor}

\usepackage{subcaption}
\usepackage{amsmath}
\usepackage{amssymb}
\usepackage{mathtools}  

\usepackage[style=numeric, sorting=none]{biblatex}
\title{Duality in Optimization}
\subtitle{A Visual Tutorial}
\author{Mohini Bariya, Rhiza Research}
\date{\today}

\begin{document}
\maketitle

\begin{center}
{\bfseries\abstractname}
\end{center}
\begin{center}
\begin{minipage}{0.9\textwidth}
\small
Lagrangian duality is fundamental to optimization, from the development of solver techniques to the interpretation of solutions; yet its geometric intuition can be difficult to grasp through mathematical derivation alone. This tutorial explains duality and its implications in a highly visual way for a broad audience of students and practitioners. It develops the motivation for the Lagrangian, the two-player game interpretation, and the connection to Fenchel duality. Throughout, it emphasizes the geometric meaning of the solution from both the primal and dual perspectives.
\end{minipage}
\end{center}

\section{Introduction}
Lagrangian duality in optimization can feel like a mathematical slight of hand. \marginpar{\footnotesize I would love to hear your feedback so I might continue to improve this tutorial. Please write me at \href{mailto:mohini@berkeley.edu}{mohini@berkeley.edu}} Each derivation step seems simple and clear; yet the final result remains cryptic. While there are many excellent duality tutorials, there seems to be a lack of visually rich explanations that build intuition beyond equations. This tutorial tries to fill that gap. It presents duality simply and visually, aiming to be broadly accessible with little prerequisite. I hope it is helpful to you and leaves you with new or stronger intuition.

\section{A Constrained Optimization Problem}
Let us say we are trying to solve the following problem: 
\begin{gather*}
    \min_x f(x)\\
    \textnormal{s.t. } x \leq 0
\end{gather*}
This is an optimization problem where we are minimizing an objective --- $f(x)$ --- under a constraint --- in this case we assume the simple constraint $x \leq 0$. For now, let us not assume any convexity [Fig. \ref{fig:cvx}].
\marginpar{
  \captionsetup{type=figure}
  \includegraphics[width=\marginparwidth]{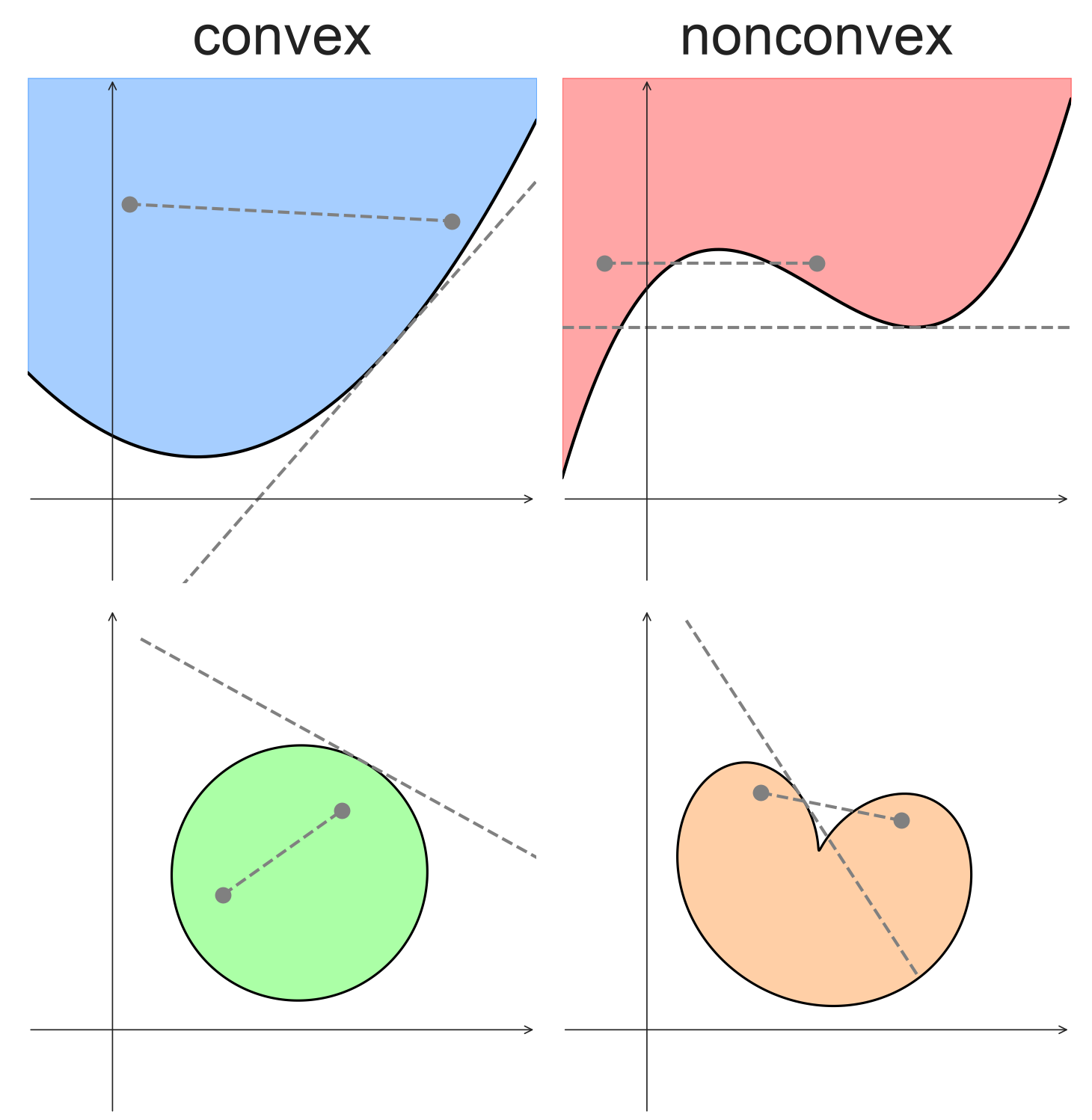}
  \caption{\footnotesize Convex functions (and sets) have nice curvature (closure) properties that makes them easy to work with. A local minimum is a global minimum; tangents lie below everywhere. This note does not discuss convexity properties extensively; see \cite{bertsekas2009convex} for a much more thorough discussion.}
  \label{fig:cvx}
}
Note that conventionally, we actually have no idea how to solve this problem. We know how to minimize a function (we could, for example, apply some gradient descent method which, depending on whether the function is convex or nonconvex, will get us the global or local minimum). However, we have \textit{no tools} to incorporate the constraint.  How can we integrate the constraint into our function in order to get back within the realm of function minimization that we are comfortable with?\\
\newline
We could define a new objective function like this: 
\begin{gather*}
    J(x) = \begin{cases}
        f(x) \textnormal{ if } x \leq 0\\
        \infty \textnormal{ otherwise}
    \end{cases}\\
    \triangleq  f(x) + I[x] \textnormal{ where } I[x] = \begin{cases}
        0 \textnormal{ if } x \leq 0\\
        \infty \textnormal{ otherwise}
    \end{cases}
\end{gather*}
Here $I[x]$ is a feasibility indicator function. Where the original optimization problem is feasible, $J(x)$ takes on the original objective value; where it is infeasible, $J(x)$ is infinite. This new function does capture our constraints, and if we minimized it, we would be solving the original problem. However, it is intractable to work with because $I[x]$ is non-differentiable. What if we replace $I[x]$ with something more tractable?
\begin{gather*}
    I(x) = -\lambda x\\
\end{gather*}
This is not perfect, but if $\lambda \leq 0$ then the penalty is at least in the right direction\footnote{many optimization presentations will flip the sign here and define $I(x) \triangleq \lambda x$ and $\lambda \geq 0$; Our less common convention makes the connection to Fenchel duals easier; as we shall see. Regardless of the chosen sign convention, the logic of $I(x)$ is the same.}: it penalizes infeasible $x$ values ($x > 0$) by contributing positively to the objective we wish to minimize, and rewards values that are less than $0$ by reducing the objective. $I(x)$ is a linear relaxation of $I[x]$ and will form a lower bound of the latter \cite{knowles2010lagrangian}, taking negatives values where $I[x]$ is 0 and finite positive values where $I[x]$ is infinite [Fig. \ref{fig:Ix}].\\ \marginpar{
  \captionsetup{type=figure}
  \includegraphics[width=\marginparwidth]{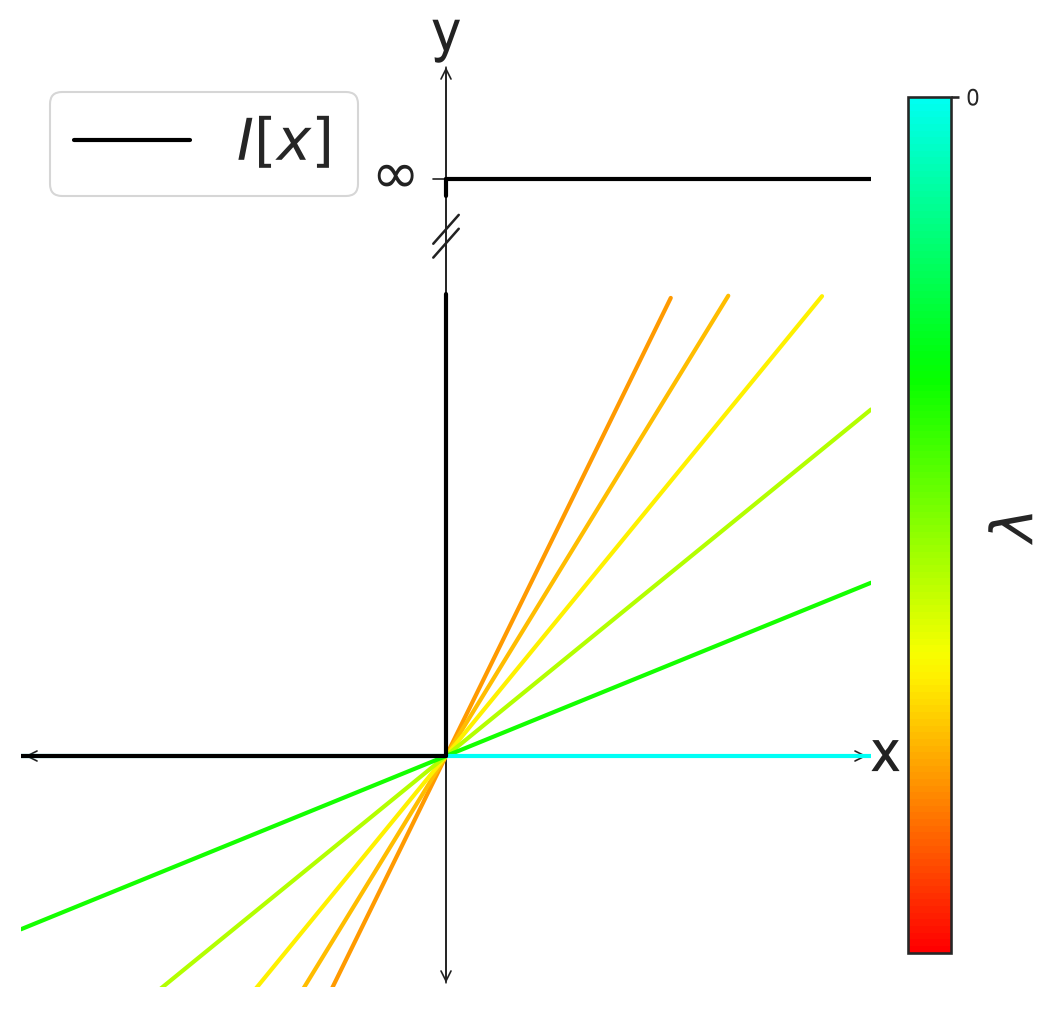}
  \caption{\footnotesize $I(x)$ for any $\lambda \leq 0$ is a lower bound on $I[x]$.}
  \label{fig:Ix}
}
\newline
This pragmatic choice of $I(x)$ gives us the \textbf{Lagrangian}: 
\begin{gather*}
    L(x, \lambda) = f(x) - \lambda x
\end{gather*}
Notice: 
\begin{gather*}
    \max_{\lambda \leq 0} L(x, \lambda) = J(x)
\end{gather*}
If $x$ is infeasible ($x > 0$), the maximizing $\lambda$ is $\lambda = -\infty$, giving $L(x, \lambda) = \infty$. If $x$ is feasible ($x \leq 0$), the maximizing $\lambda$ is $\lambda = 0$, giving $L(x, \lambda=0) = f(x)$. So, our original optimization problem can now be rewritten as: 

\[
\left.
\begin{gathered}[c]
    \min_x f(x) \\[4pt]
    \text{s.t. } x \leq 0
\end{gathered}
\quad\Longleftrightarrow\quad
\min_x \max_{\lambda \leq 0} L(x, \lambda)
\;\right\}\;\text{primal problem}
\]
We refer to this as the \textbf{primal problem}. What does the min-max order mean? We are choosing an $x$ such that the value of $L$ is minimized over all choices of $\lambda$. In other words, choose an $x$ such that the worst-case outcome over all $\lambda$ is minimized. In a game-theory framing, player 1 starts by picking an $x$ such that player 2's choice of $\lambda$ will do the least damage (maximization) to the objective. \\
\newline
This is still a hard problem to solve because we are trying to minimize robustly against a set of maximizers that effectively resemble the discontinuous $I[u]$. What if we switch the order?
\[
\left.
\max_{\lambda \geq 0} \min_x L(x, \lambda)
\;\right\}\;\text{dual problem}
\]
Switching the min-max for a max-min will give us a lower bound on the optimization: “[The] minimums of maximums are always greater than maximums of minimums” \footnote{ \footnotesize See Ben Recht's duality blog for this and more \cite{recht2024duality}.}. This is called the \textbf{dual problem} and turns out to be more tractable to work with, as we will see in the following sections. 
\subsection{Visualizing the Lagrangian \& the dual}
\begin{figure}[b!]
  \centering
  \includegraphics[width=0.8\textwidth]{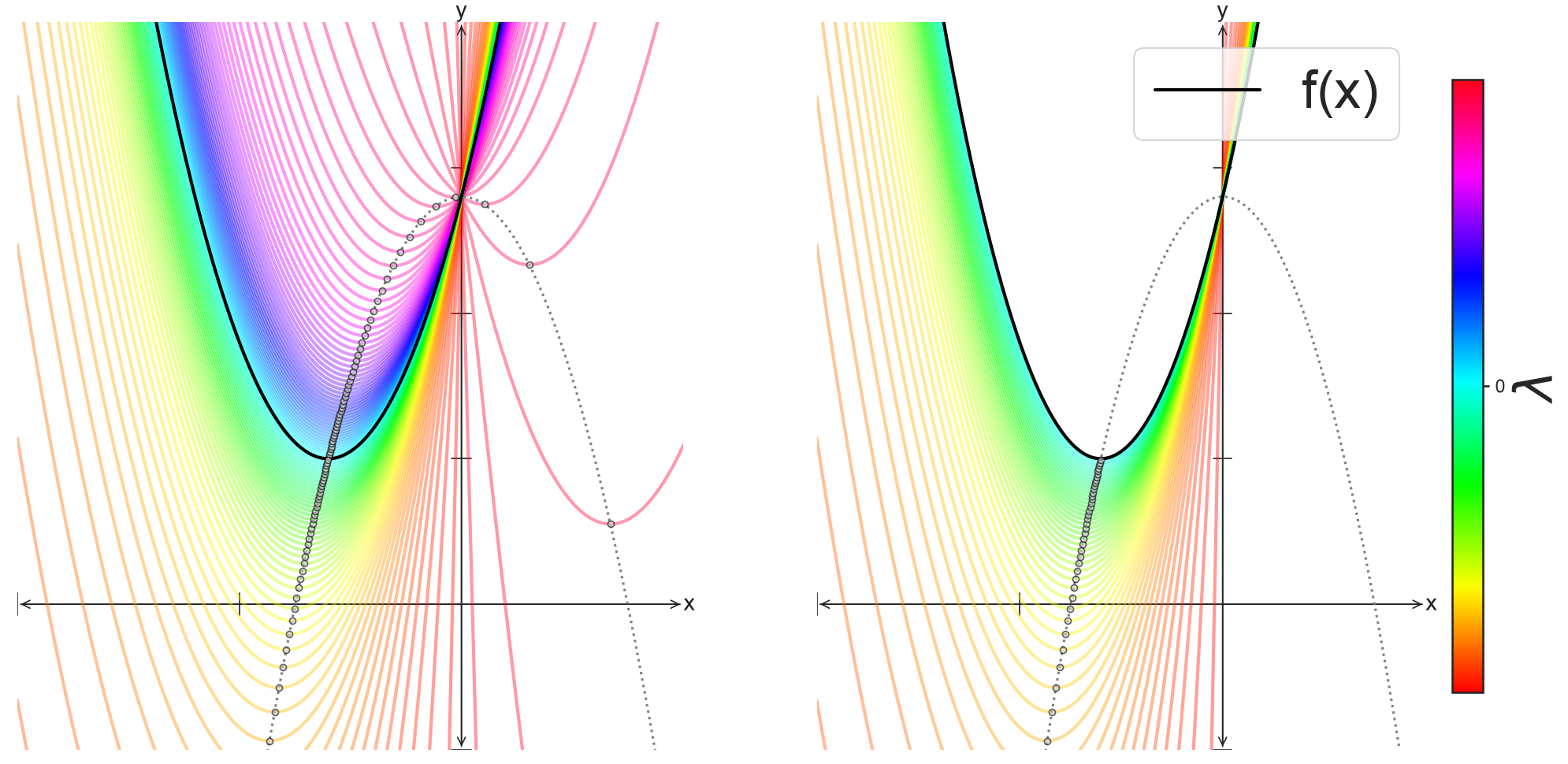}
  \caption{\footnotesize The Lagrangian of $f(x) \textnormal{ st } x \leq 0$ - across all values of $\lambda$ (left) and only for $\lambda \leq 0$ (right) .}
  \label{fig:Lagrangian}
\end{figure}

We can build intuition on the Lagrangian by visualizing it. We take the objective function and the constraint and combine them into a single function, through a penalty weight $\lambda$. The result is a slightly different function for each value of $\lambda$ [Fig. \ref{fig:fminusl}]. \marginpar{
  \captionsetup{type=figure}
  \includegraphics[width=\marginparwidth]{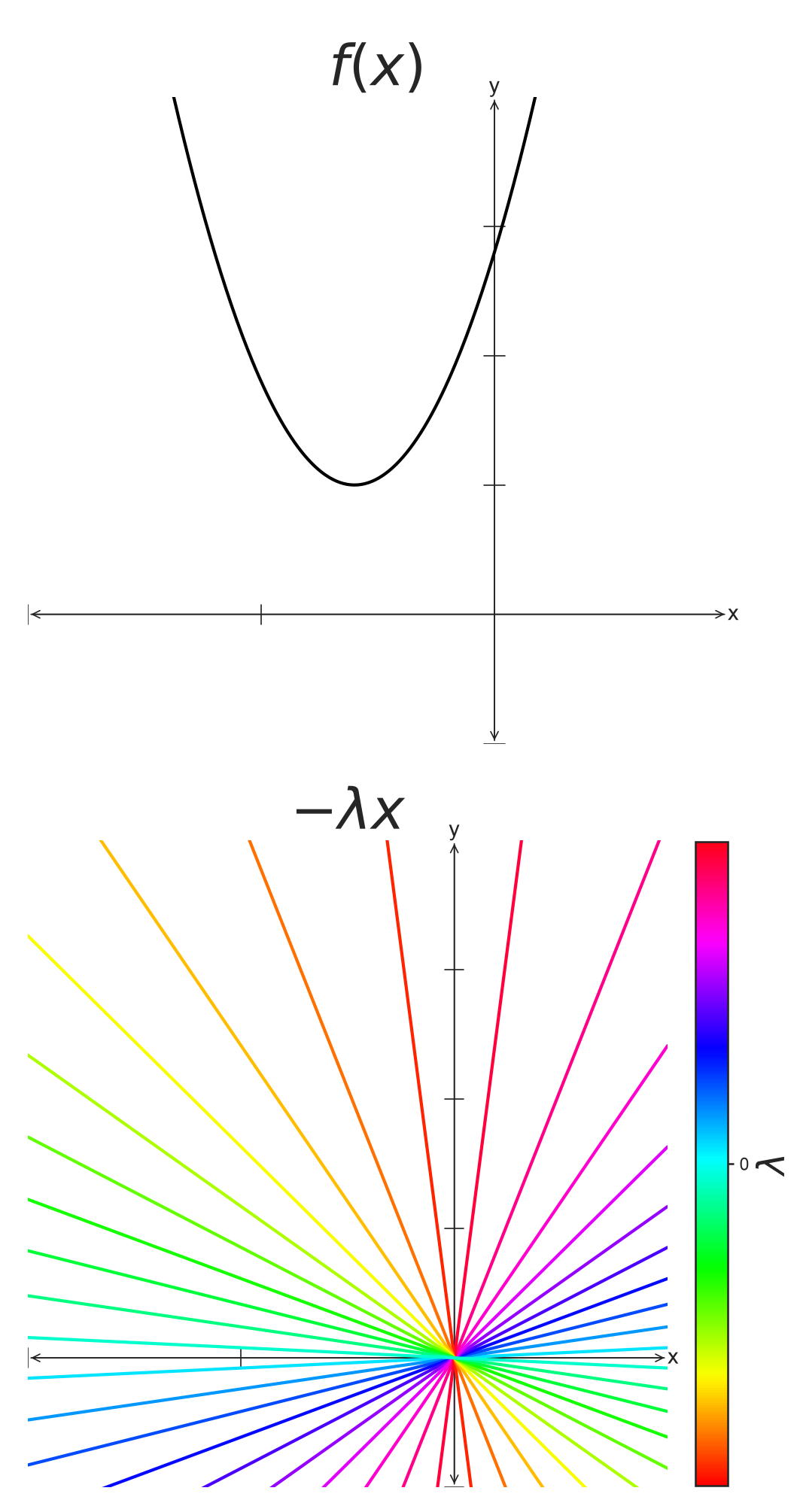}
  \caption{\footnotesize Objective \& constraint are combined with penalty weight $\lambda$.}
  \label{fig:fminusl}
} Collecting these across $\lambda$, we see that the Lagrangian consists of a set of curves related to the original objective $f(x)$ [Fig. \ref{fig:Lagrangian}]. 
$\min_x L(x, \lambda)$ picks up a single point on each curve; collectively these minima across all $\lambda$ values draw out a new curve. We can see from Fig. \ref{fig:Lagrangian} that the resulting curve is concave (this is a powerful result that Section \ref{sec:dualisnicer} shows  holds in general). It intersects the original objective $f(x)$ at the global optimum --- for $\lambda = 0$ --- and at the boundary of the feasible region --- in this case $ x = 0$ \footnote{\footnotesize Thanks to Tim Wheeler for creating very similar visuals to the Figures in this section \cite{wheeler2022convexduality}}.

\subsubsection{Primal and dual as two-player games}
\begin{figure}[h!]
  \centering
  \includegraphics[width=\linewidth]{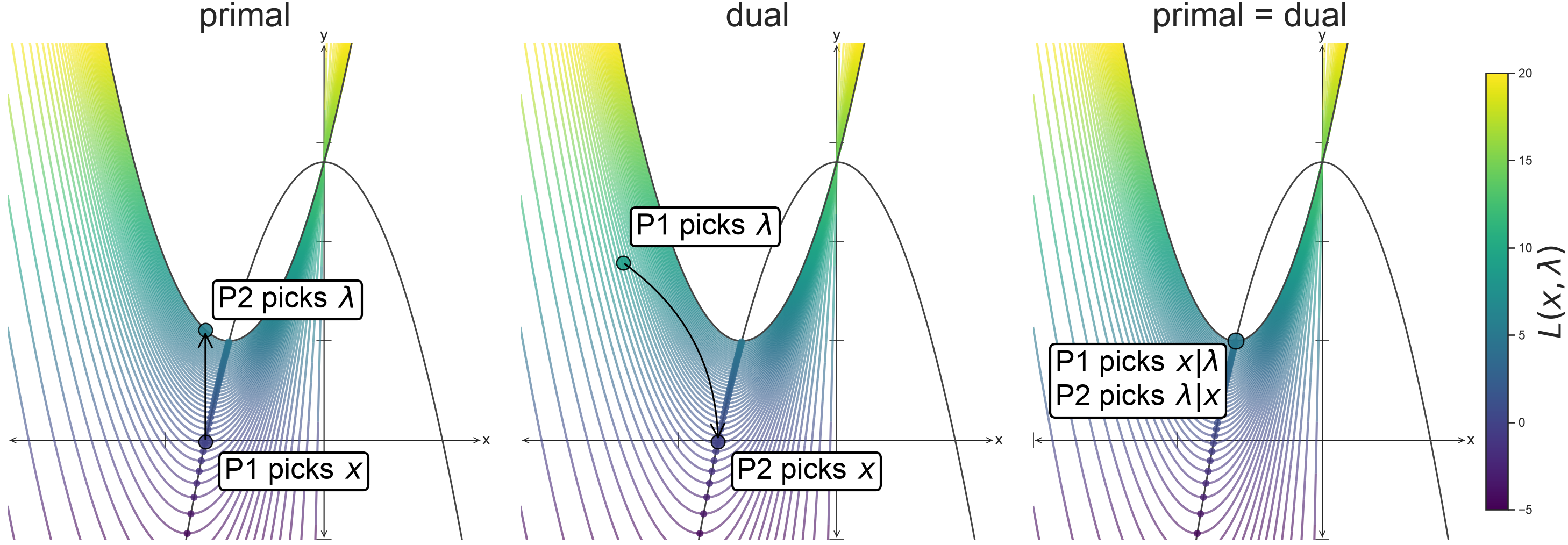}
  \caption{\footnotesize Visualizing the Lagrangian as a two player game. Curves are colored by value of $L(x, \lambda)$.}
  \label{fig:game}
\end{figure}

This visualization helps us understand the interpretation of the primal and dual problems as two player games [Fig. \ref{fig:game}]. At left is the primal view, $\min_x \max_\lambda L(x, \lambda)$: Player 1 picks $x$ first trying to minimize the objective; Player 2 then moves across the $\lambda$-curves at that $x$, pushing to the highest possible point — which lands on $f(x)$. At center is the dual view, $\max_\lambda \min_x L(x, \lambda)$: Player 1 picks $\lambda$ first trying to maximize the objective; Player 2 then moves along the curve traced by that $\lambda$, pushing to the lowest possible point. In both the primal and dual cases, Player 1's choice is suboptimal, giving Player 2 more room to counteract than necessary. At right is the choice that is optimal regardless of whether Player 1 is trying to minimize over $x$ or maximize over $\lambda$ --- it is better at minimizing the objective than the $x$ chosen at left and better at maximizing the objective than the $\lambda$ chosen at center. This is exactly the solution of the original optimization problem.

\subsubsection{The dual problem is nicer}\label{sec:dualisnicer}
This is all well and good but why bother with the dual problem at all? The dual problem is nicer than the primal.
\begin{gather*}
    \textnormal{ let } g(\lambda) = \min_x \bigg[f(x) - \lambda x\bigg]
\end{gather*}
This is a pointwise minimum (over $x$) across a set of functions that are affine in $\lambda$; therefore $g(\lambda)$ must be concave, regardless of the shape of the original $f(x)$. Since the dual problem is $\max_{\lambda} g(\lambda)$, a maximization of a concave function, it is always tractable.\\ 
\newline
Fig. \ref{fig:xlambda} visualizes the Lagrangian versus $x$ and versus $\lambda$ - we can clearly see how the pointwise minimum over a set of affine functions results in a concave function. You can work out that the minimum points along each $\lambda$ curve at top correspond to the pointwise minimum across affine functions of $\lambda$ at bottom -- these are just two different views of $L(x, \lambda) = \min_x f(x) - \lambda x$.
\marginpar{
  \captionsetup{type=figure}
  \includegraphics[width=\marginparwidth]{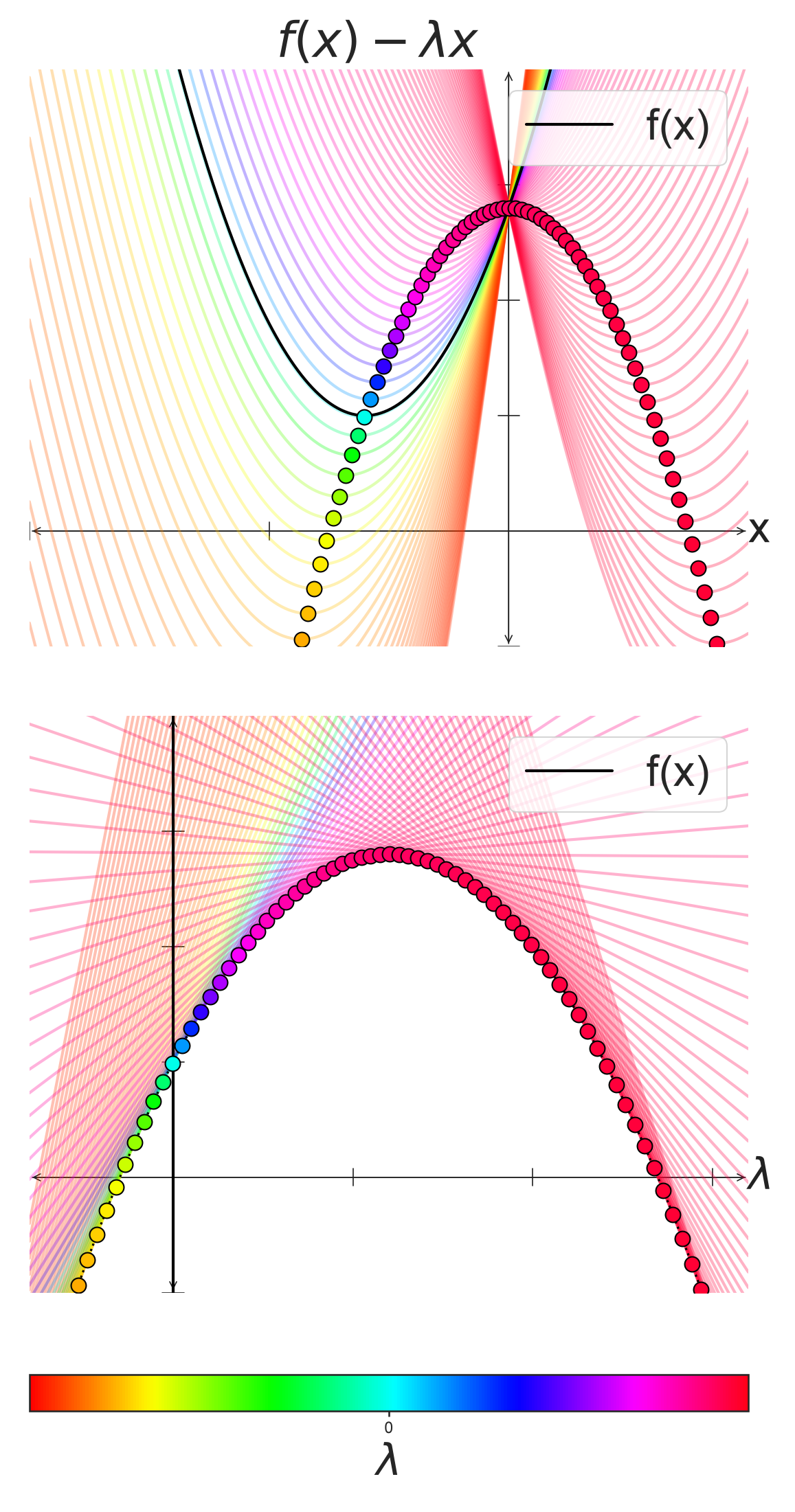}
  \caption{\footnotesize Two views of $L(x, \lambda)$: versus $x$ (top) and versus $\lambda$ (bottom). The colored points correspond and are exactly $\min_x f(x)-\lambda x$ for different values of $\lambda$. $f(x)$ in the top plot corresponds to the y-axis in the bottom plot where $\lambda = 0$.}
  \label{fig:xlambda}
}
Fig. \ref{fig:nonconvexf} visualizes $f(x)$ and $g(\lambda)$ when $f(x)$ is obviously nonconvex; we can see that nevertheless $g(\lambda)$ is concave and tractable. \\
\marginpar{
  \captionsetup{type=figure}
  \includegraphics[width=\marginparwidth]{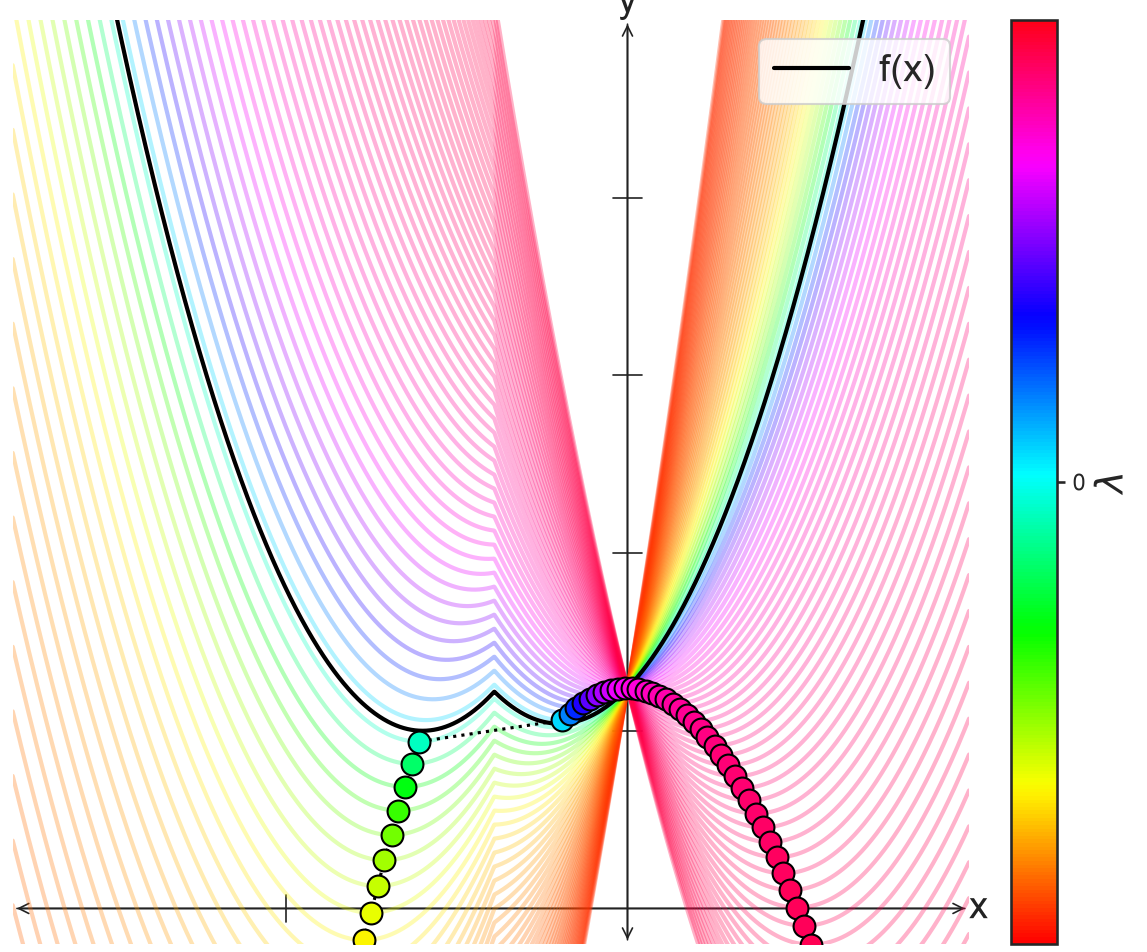}
  \caption{\footnotesize $g(\lambda)$ remains concave even when $f(x)$ is nonconvex.}
  \label{fig:nonconvexf}
}
\newline
In the next section we will see that the dual problem has another interpretation which makes intuitive many of the interesting results we obtain in solving the dual problem.
\section{The Fenchel Dual}
\begin{figure*}[tbp]
  \centering
  \includegraphics[width=\linewidth]{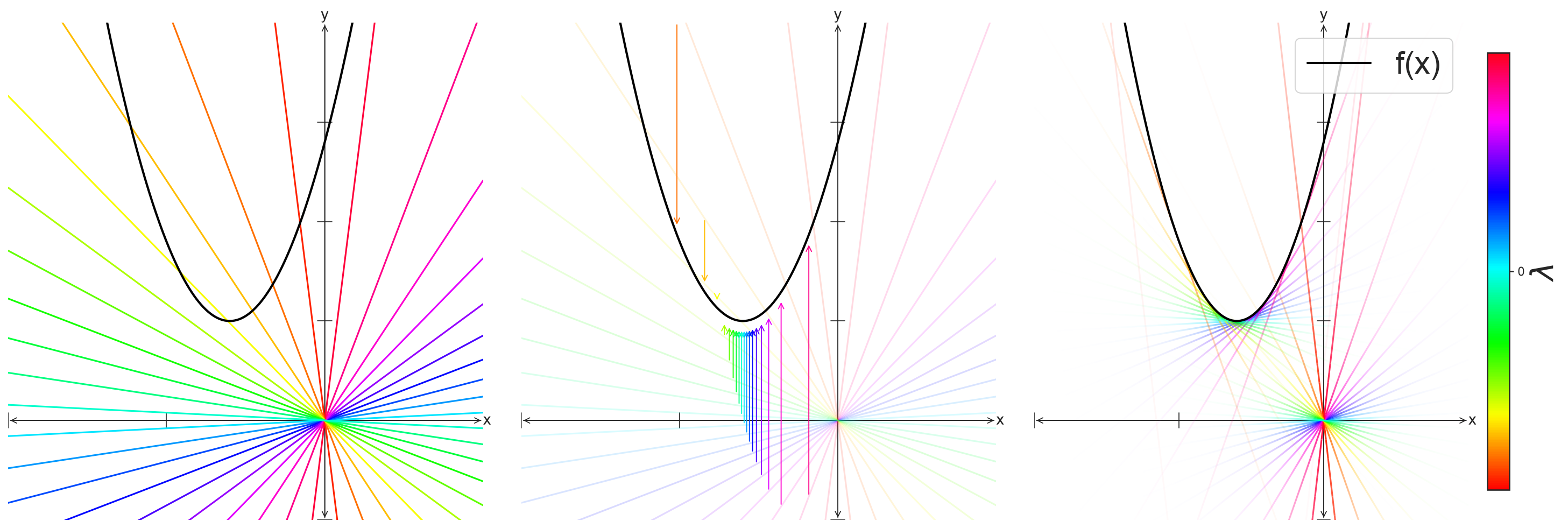}
  \caption{The Fenchel dual parametrizes the supporting hyperplanes of $f(x)$; for every line $\lambda x$, the Fenchel dual equation determines how much it must be vertically shifted to become a supporting hyperplane of $f(x)$.}
  \label{fig:supporting-hyperplanes}
\end{figure*}

Given a function $f(x)$, the \textbf{Fenchel conjugate} or dual is defined as: 
\begin{gather*}
    f^*(\lambda) = \max_x \lambda x - f(x)
\end{gather*}
That looks familiar! It is in fact $-g(\lambda)$, the inner part of our dual problem: 
\begin{gather*}
    - g(\lambda) = -\bigg[\min_x f(x) - \lambda x\bigg] = \max_x \lambda x - f(x) = f^*(\lambda)
\end{gather*}
where the negation flips the $\min$ into a $\max$. In mathematics, a dual allows us to view a space or a function from another perspective that preserves the same information, but could be more workable. What is this new perspective in the case of the Fenchel conjugate and --- by extension --- the dual optimization problem?\\
\newline
In fact, the Fenchel conjugate of a function parametrizes the set of hyperplanes that support the function (in a single dimension, these hyperplanes are tangent lines). For a truly convex function, the supporting hyperplanes fully describe the function. \\
\newline
Fig. \ref{fig:supporting-hyperplanes} illustrates why the equation defining the Fenchel conjugate indeed parameterizes the set of supporting hyperplanes. For every $\lambda$, there is a supporting hyperplane to convex function $f(x)$ of slope $\lambda$. To find this hyperplane, we need to move the line $\lambda x$ up or down by a certain amount so that it just grazes the function $f(x)$. This amount is given by: 
\begin{gather*}
    \Delta = \max_x \lambda x - f(x)
\end{gather*}
Lines that pass above $f(x)$ will be shifted down by their largest gap to $f(x)$; lines passing below will be shifted up by their smallest gap to $f(x)$. The Fenchel dual is exactly given by the above equation, and --- over all $\lambda$ --- parametrizes the full set of supporting hyperplanes that define $f(x)$. 
\begin{gather*}
f^*(\lambda) = \max_x \lambda x - f(x) = \Delta_{\lambda}
\end{gather*}
$f^*(\lambda)$ is convex in $\lambda$ [Fig. \ref{fig:dualvslambda}]. But we already know this! After all, $g(\lambda) = -f^*(\lambda)$ and we have seen from the above discussion that $g(\lambda)$ is concave.
\newline
\marginpar{
  \captionsetup{type=figure}
  \includegraphics[width=\marginparwidth]{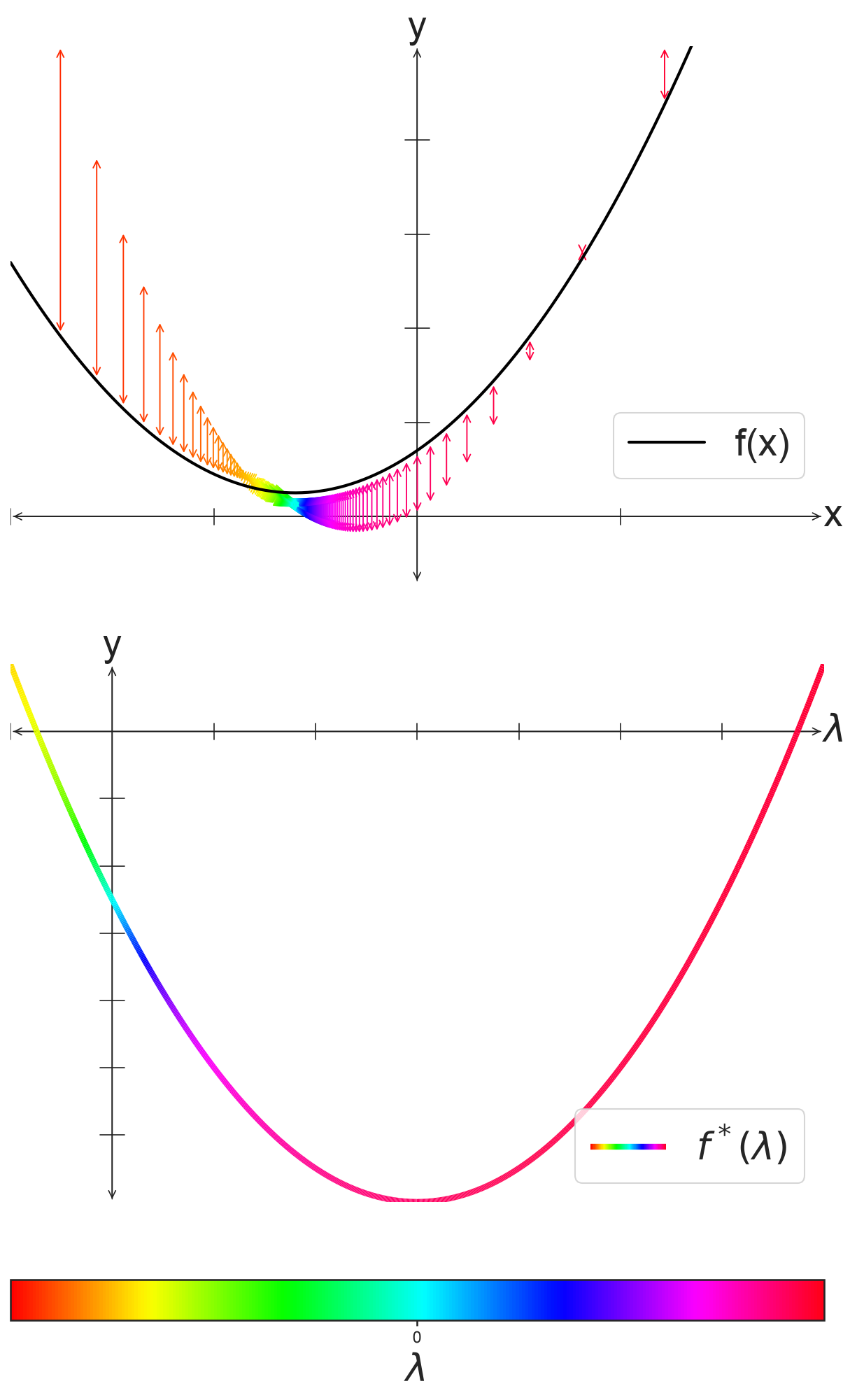}
  \caption{\footnotesize The Fenchel dual parametrizes the supporting hyperplanes of $f(x)$ by picking out the vertical shifts of each hyperplane with slope $\lambda$. The result when plotted versus $\lambda$ is a convex function.}
  \label{fig:dualvslambda}
}
\subsection{Optimal point}
\begin{figure}[h!]
  \centering
  \includegraphics[width=\linewidth]{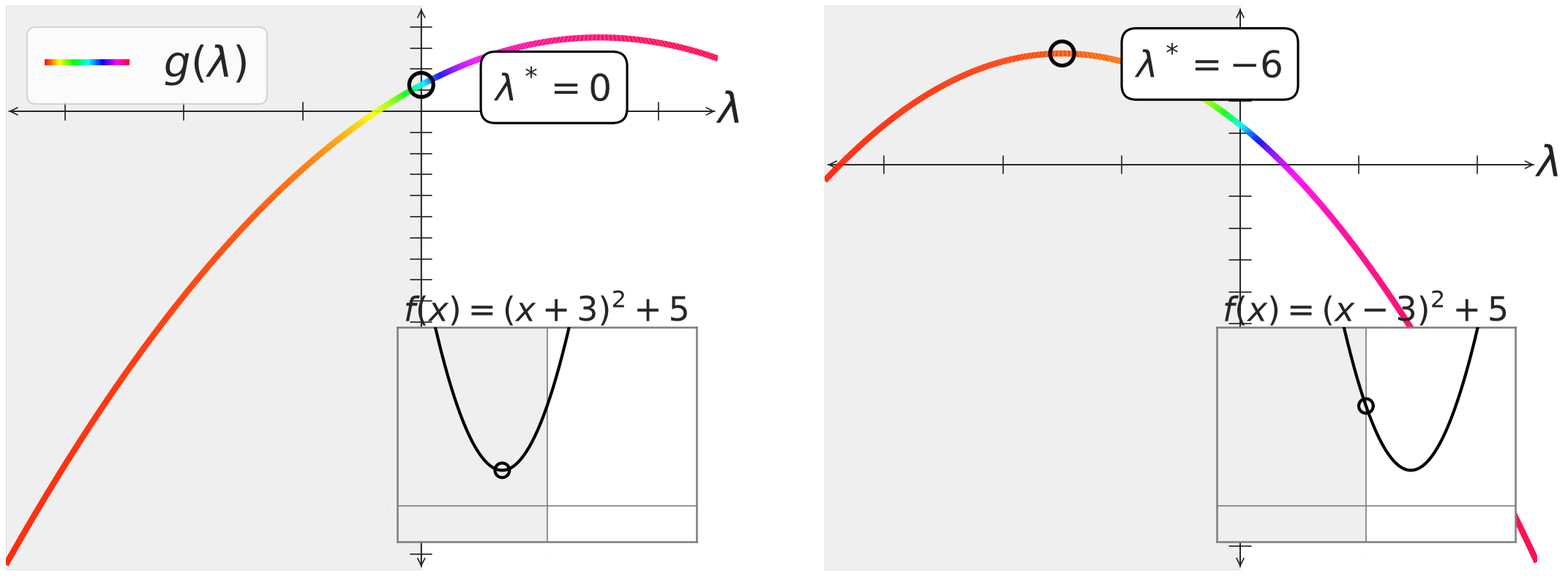}
  \caption{\footnotesize Visualizing the dual solution $\max_{\lambda \leq 0} g(\lambda)$ when the primal solution is unconstrained (left) and constrained (right).}
  \label{fig:cvsuc}
\end{figure}

What can we say about the optimal solution when viewed from this dual perspective of supporting hyperplanes? The solution of the dual problem is $\lambda^* = \arg\max_{\lambda \leq 0} g(\lambda)$, where $g(\lambda) = -f^*(\lambda)$. Since $g(\lambda)$ is concave, this maximum will occur either at the unconstrained maximizer of $g$ \textit{or} at the boundary of the constraint ($\lambda = 0$). We can go through both cases.\\
\newline
\textbf{Unconstrained solution.} Let us first locate the
unconstrained maximizer, being careful to remember that $x$ is a function of $\lambda$. Defining $x^*(\lambda) = \arg\min_x\big[f(x) - \lambda x\big]$, we have:
\begin{gather*}
    \frac{d}{d\lambda}g(\lambda)
    = \frac{d}{d\lambda}\Big[f(x^*(\lambda)) - \lambda x^*(\lambda)\Big]
    = -x^*(\lambda)^{\textnormal{[see Appendix \ref{app:envelope}]}}\\
    \textnormal{maximum occurs at: } \frac{d}{d\lambda}g(\lambda) = 0\\
    \implies x^*(\hat\lambda) = 0
\end{gather*}
Let us call this unconstrained maximizer $\hat\lambda$. This is exactly what we see in Fig.~\ref{fig:xlambda}: the peak of $g(\lambda)$ corresponds to the point $x=0$, and $\hat\lambda$ is the derivative of $f(x)$ at $x=0$.\footnote{That the maximum of $g(\lambda)$ will occur at the constraint boundary can be shown for a more general constraint as well. If our original problem is $\min_x f(x) \textnormal{ s.t. } c(x) \leq 0$, the same argument
gives $\frac{d}{d\lambda}g(\lambda) = -c(x^*(\lambda))$; thus at the unconstrained maximizer, $c(x^*(\hat\lambda)) = 0$.}\\
If $\hat\lambda$ lies within our feasible region $\lambda \leq 0$, it is the solution of the dual problem: $\lambda^* = \hat\lambda$. If it doesn't...\\
\newline
\textbf{Constrained solution.} We can show that the
maximizing $\lambda$ must be $\lambda^* = 0$, as follows.
$\hat\lambda$ outside the feasible region means $\hat\lambda > 0$,
so the peak of $g(\lambda)$ lies to the right of the origin. Since
$g(\lambda)$ is concave, it has only one peak, so it must be
increasing everywhere to the left of $\hat\lambda$ -- in particular,
increasing throughout the whole feasible region $\lambda \leq 0$.
The maximum over this feasible region therefore occurs at its
right-hand boundary, $\lambda^* = 0$, which corresponds to
$g(0) = \min_x f(x)$, the unconstrained optimum of the original
problem.\\
\newline
Here we have another echo of "duality". The constrained solution of the dual corresponds to the unconstrained solution of the primal, while the unconstrained solution of the dual corresponds to the constrained solution of the primal.
\begin{aside}
    Succinctly: the solution of our dual problem --- let us denote it $(x^*, \lambda^*)$ --- will either lie at the boundary of the original constraint (in our case, $x^* = 0$) with corresponding $\lambda^* < 0$, \textit{or} it will lie exactly at $x^* = \min_x f(x)$ with corresponding $\lambda^* = 0$. $\lambda^*$ is also the gradient of $f(x)$ at $x^*$. 
\end{aside}Fig. \ref{fig:cvsuc} visualizes this result for the two cases: where the primal solution is constrained versus unconstrained. We can also see that the gradient of $f(x)$ at the solution is $\lambda^*$. That $\lambda^* = 0$ at the unconstrained solution and $\lambda^* \neq 0$ at the constrained solution is called "complementary slackness". 

\subsection{$\lambda^*$ is the "shadow price"}
Many times in applied optimization --- particularly in economics applications --- you will see $\lambda^*$ referred to as the shadow price. The hyperplane view makes the reason quite clear. Typically in complex, real-world optimization problems, the solution will not lie at the unconstrained point, but at some boundary point where one or more constraints are active. This point will have $(x^*, 
\lambda^* \neq 0)$, and $\lambda^*$ is the local gradient of the objective $f(x)$ at that boundary. Therefore, $\lambda^*$ tells us how much we could shift the objective value if only we were allowed to budge $x^*$, i.e. to loosen the binding constraint.

\printbibliography

\appendix
\section{Derivation of $\frac{d}{d\lambda}g(\lambda) = -x^*(\lambda)$}
\label{app:envelope}

Let $x^*(\lambda) = \arg\min_x \big[f(x) - \lambda x\big]$, so that
$g(\lambda) = f(x^*(\lambda)) - \lambda x^*(\lambda)$. Differentiating
with respect to $\lambda$:
\begin{gather*}
    \frac{d}{d\lambda}g(\lambda)
    = f'(x^*(\lambda))\,\frac{dx^*}{d\lambda} - x^*(\lambda) - \lambda\,\frac{dx^*}{d\lambda} \\
    = \Big[f'(x^*(\lambda)) - \lambda\Big]\frac{dx^*}{d\lambda} - x^*(\lambda)
\end{gather*}
By the first-order condition defining $x^*(\lambda)$,
\begin{gather*}
    f'(x^*(\lambda)) - \lambda = 0
\end{gather*}
so the bracketed term vanishes, giving
\begin{gather*}
    \frac{d}{d\lambda}g(\lambda) = -x^*(\lambda)
\end{gather*}
\end{document}